\documentclass[10pt,twoside]{article}
\usepackage{amsmath,amsbsy,amsfonts}
\newtheorem{lemma}{Lemma}[section]

\newtheorem{theorem}{Theorem}[section]

\newtheorem{coro}{Corollary}[section]

\newcommand{\fim}{\hfill\rule{2mm}{2mm}}

\newcommand{\R}{\mbox{${\rm{I\!R}}$}}

\let\Section=\section
\def\section{\setcounter{equation}{0}\Section}
\begin{document}
\title{Multivalued Elliptic Equation with exponential critical growth in $\mathbb{R}^2$ \thanks{Partially supported by Procad and Casadinho CNPQ } }

\author{Claudianor O. Alves\thanks{Partially supported by  CNPq - Grant 304036/2013-7 }\,\, and \,\,Jefferson A. Santos   }

\date{}

\pretolerance10000

\maketitle

\begin{abstract}
In this work we study the existence of nontrivial solution for the
following class of multivalued elliptic problems
$$
-\Delta u+V(x)u-\epsilon h(x)\in \partial_t F(x,u) \quad \text{in} \quad \mathbb{R}^2, \eqno{(P)}
$$
where $\epsilon>0$, $V$ is a continuous function verifying some conditions, $h \in (H^{1}(\mathbb{R}^{2}))^{*}$ and 
$\partial_t F(x,u)$ is a generalized gradient of $F(x,t)$ with
respect to $t$ and $F(x,t)=\int_{0}^{t}f(x,s)\,ds$. Assuming that $f$ has an exponential critical growth and a discontinuity point, we have applied Variational Methods for locally Lipschitz functional to get two solutions for $(P)$ when $\epsilon$ is small enough. \\



\end{abstract}

\section{Introduction}

The interest in the study of nonlinear partial differential equations with discontinuous
nonlinearities has increased because many free boundary problems arising
in mathematical physics may be stated in this from. Among these problems, we
have the obstacle problem, the seepage surface problem, and the Elenbaas equation;
see for example \cite{Chang, Chang2, Chang3}.

Among the typical examples, we have chosen the model for the heat conductivity
in electrical media. This model has a discontinuity in its constitutive laws. In fact,
considering a domain $\Omega \subset \mathbb{R}^{2}$ (which in particular could be taken as being the whole space
$\mathbb{R}^{2}$, see  \cite{ACD}) with electrical media, the thermal and electrical conductivity are denoted by
$K(x, t)$ and $\sigma(x, t)$, respectively. Here $x$ is in $\Omega$ and $t$ represents the temperature.
Since we are considering an electrical media, the function $\sigma$ may have discontinuities
in $t$, and the distribution of the temperature is unknown. The differential equation
describing this distribution is
$$
-\sum_{i=1}^{n}\frac{\partial}{\partial x_i}\left((K(x,u(x))\frac{\partial u(x)}{\partial x_i}\right)=\sigma(x,u(x)).
$$

Note that this equation is related to a free boundary problem in which the jump
surface of the electrical conductivity is unknown. We describe this surface as being
the set
\begin{equation} \label{gamma}
\Gamma_\alpha(u)=\{x \in \Omega\,:\,u(x)=\alpha, \quad \sigma \, \mbox{is discontinuous at} \, \alpha \}.
\end{equation}

When the thermal conductivity $K$ is constant, $\Omega=\mathbb{R}^{2}$ and the electrical conductivity $\sigma$ is given by 
$$
\sigma(x,t)=H(t-a)f(t) +\epsilon h(x)-V(x)t,
$$
with $f$ having an exponential critical growth, the model  becomes 
\begin{equation}\label{Peps1}
-\Delta u + V(x)u= H(u-a)f(u) +\epsilon h(x) \quad \mbox{in} \quad \mathbb{R}^{2}.
\end{equation}
Here $H$ is the Heaviside function, that is, 
$$
H(t)=
\left\{
\begin{array}{l}
1, \quad t>0,\\
0, \quad t \leq 0,
\end{array}
\right.
$$
$\epsilon$ is a positive parameter and $h$ is a measurable function defined in $\mathbb{R}^2$. 
Note that in this model the jump surface of the solution (\ref{gamma}) is represented by
the set
\begin{equation} \label{gamma2}
\Gamma_a(u)=\{x \in \mathbb{R}^2 \,:\,u(x)=a \}.
\end{equation}

Related to Problem (\ref{Peps1}) for the special case of $a = 0$, i.e., without jump discontinuities,
we cite the works by de Freitas \cite{freitas} and do \'O, de Medeiros and Severo \cite{JEU1, JEU2}. 

A rich literature is available by now on problems with discontinuous nonlinearities,
and we refer the reader to  Alves, Bertone and Gon\c calves \cite{ABG}, Alves and Bertone \cite{AB}, Alves, Gon\c calves and Santos \cite{Abrantes},  , Ambrosetti and Turner \cite{turner}, Ambrosetti, Calahorrano and Dobarro \cite{ACD}, Badiale and Tarantelo \cite{Badiale}, Carl, Le and Motreanu \cite{carl}, Clarke \cite{clarke}, Chang \cite{Chang}, Carl and Dietrich \cite{CD1}, Carl and S. Heikkila \cite{CD2,CD3}, Cerami \cite{cerami}, de Souza, de Medeiros and Severo \cite{MES1,MES2}, Hu, Kourogenis and Papageorgiou \cite{Hu}, Montreanu and Vargas \cite{MV}, Radulescu \cite{Radulescu} and their references. Several techniques have been developed or applied in their study, such as variational methods for nondifferentiable functionals, lower and upper solutions, global branching, fixed point theorem, and the theory of multivalued mappings.

After a  bibliography review, we did not find any paper involving existence of solution for a class of elliptic problem with discontinuous nonlinearity and exponential critical growth via  variational methods for nondifferentiable functional. Motivated by this fact, in this paper we employ variational techniques to study existence and multiplicity of nonnegative solutions for a large class of multivalued elliptic equations, which includes the equation (\ref{Peps1}).  More precisely, we will study the multivalued elliptic equation
$$
-\Delta u+V(x)u-\epsilon h(x)\in \partial_t F(x,u) \quad \text{in} \quad \mathbb{R}^2, \eqno{(P)}
$$
where $\epsilon >0$ is a positive parameter, $V$ is a continuous function verifying some technical conditions, $h \in (H^{1}(\mathbb{R}^{2}))^{*}$ and $F(x,t)$ is the primitive of a function $f(x,t)$, which has an exponential critical growth and a discontinuity point, for more details see Section 2.

In  $\mathbb{R}^2$, to apply variational methods, the natural growth restriction on the function $f$  is given by the inequality
of Trudinger and Moser \cite{M,T}. More precisely, we say that a
function $f$ has an exponential critical growth if there is $\alpha_0 >0$ such that
$$ \lim_{|s| \to \infty} \frac{|f(s)|}{e^{\alpha s^{2}}}=0
\,\,\, \forall\, \alpha > \alpha_{0}\quad \mbox{and} \quad
\lim_{|s| \to \infty} \frac{|f(s)|}{e^{\alpha s^{2}}}=+ \infty
\,\,\, \forall\, \alpha < \alpha_{0}.
$$
We would like to mention that problems involving exponential critical growth, with $f$ being a continuous functions, have received a special attention at last years, see for example, \cite{1, AdoOM,AP, Cao, DMR, DdoOR, doORuf,  5} and their references. Here, since we intend to get a solution for the differential inclusion $(P)$, we assume that there exists $\alpha_0>0$ such that 
\begin{description}
\item[$(f_0)$]$\displaystyle \limsup_{t\to +\infty}\frac{\max\left\{|\xi|;\xi\in \partial_{t}F(x,t)\right\}}{e^{\alpha_0|t|^{2}}}<+\infty $  uniformly in $x \in \mathbb{R}^2.$
\end{description}
Moreover, assuming a condition at origin like 
\begin{description}
\item[$(f_1)$]$\displaystyle \limsup_{t\to 0}\frac{2\max\left\{|\xi|;\xi\in \partial_{t}F(x,t)\right\}}{|t|}<+\infty $  uniformly in $x \in \mathbb{R}^2.$
\end{description}
it is easy to check that the functional $\Psi:H^{1}(\mathbb{R}^{2}) \to \mathbb{R}$ given by
$$
\Psi(u)=\int_{\mathbb{R}^{2}}F(x,u)\,dx
$$
is well defined, for more details see Section 4. However, to apply variational methods is better to consider the functional $\Psi$ in a more appropriated domain, that is, $\Psi:L^{\Phi}(\mathbb{R}^{2}) \to \mathbb{R}$, for $\Phi(t)=e^{|t|^{2}}-1$. But,  once  $\Phi$ does not satisfy the $\Delta_2$-condition, we cannot guarantee that given $J \in (L^{\Phi}(\mathbb{R}^{2}))^*$, then 
$$ 
J(u)=\int_{\mathbb{R}^{2}}vu\,dx, \quad \forall u \in L^{\Phi}(\mathbb{R}^{2}),
$$
for some $v:\mathbb{R}^{2} \to \mathbb{R}$ mensurable function. For the familiar readers with the study of the differential inclusions, they will observe that the above remark is bad to apply variational methods, because in general for this type of equations we need to prove that the inclusion below holds
$$
\partial \Psi(u) \subset \partial_t F(x,u)=[\underline{f}(x,u(x)),\overline{f}(x,u(x))]  \text{ a.e. in } \mathbb{R}^{2},
$$
where
$$
\underline{f}(x,t)=\displaystyle \lim_{r\downarrow 0}\mbox{ess
inf}\left\{ f(x,s);|s-t|<r\right\} 
$$
and
$$
\overline{f}(x,t)=\lim_{r\downarrow 0}\mbox{ess sup}\left\{
f(x,s);|s-t|<r\right\} .
$$
 
In Section 4, we analyze this question. In fact, we show that it is enough to consider $\Psi:E_{\Phi}(\mathbb{R}^{2}) \to \mathbb{R}$ where  
$$
E_\Phi(\Omega)=\overline{C^{\infty}_0(\mathbb{R}^{2})}^{\|\,\,\|_\Phi}.
$$

Before to state our main result, we must mention our conditions on $V,h$ and $f$, which are the following:  
\begin{description}
\item[$(h_0)$] $h \in (H^{1}(\mathbb{R}^{2})))^{*}$ and $0< \displaystyle \int_{\mathbb{R}^{2}}h\,dx<+\infty$.
\item[$(V_1)$] $V$ is continuous and $V(x)\geq V_0>0$, $\forall x\in\mathbb{R}^2$,
\item[$(V_2)$] $\frac{1}{V}\in L^1(\mathbb{R}^2)$.
\end{description}

\begin{description}

\item[($f_2$)] There is $t_0\geq0$ such that 
$$
f(x,t)=0  \quad \text{for } \quad  t< t_0  \text{ and } \, \forall x \in \mathbb{R}^{2}
$$
and
$$
f(x,t)>0 \quad \mbox{for} \quad  t>t_0 \text{ and } \, \forall x \in \mathbb{R}^{2}.
$$
\item[$(f_3)$]$\displaystyle \limsup_{t\to 0}\frac{2\max\left\{|\xi|;\xi\in \partial_{t}F(x,t)\right\}}{|t|}<\lambda_1 $  uniformly in $x \in \mathbb{R}^2,$ where
$$
\lambda_1=\inf_{u \in E \setminus \{0\}}\frac{\int_{\mathbb{R}^{2}}(|\nabla u|^{2}+V(x)|u|^{2})\,dx}{\int_{\mathbb{R}^{2}}|u|^{2}\,dx}
$$
and
$$
E:=\left\{u\in H^1(\mathbb{R}^2)\,:\, \int_{\mathbb{R}^2}V(x)u^2dx<+\infty\right\}.
$$
\item[$(f_4)$] There is a compact set $K\subset \mathbb{R}^2$ and constants $c_3,c_4>0$ and $\nu>2$, such that 
$$
F(x,t)\geq c_3 t^\nu- c_4,  \quad \mbox{for}  \quad t\geq 0 \text{ and } \quad \forall x\in K.
$$
\item[$(f_5)$] There is $\tau>2$ verifying 
$$
0\leq\tau F(x,t)\leq \underline{f}(x,t)t,  \quad \mbox{for}  \quad t\geq t_0 \quad \mbox{and} \quad  \forall x\in\mathbb{R}^2.
$$
\item[$(f_6)$] There are $p>2$ and $\mu>0$ such that 
$$
F(x,t)\geq \mu (t-t_0)^p, \quad \mbox{for}  \quad t\geq t_0 \text{ and } \forall x\in\mathbb{R}^2.
$$
\end{description}
Here, we would like to point out that the function
\begin{equation} \label{funcao}
f(t)=2H(t-a)te^{t^{2}}, \quad \forall t \in \mathbb{R},
\end{equation}
verifies $(f_1)-(f_6)$.

Now, we are able to state our main result

\begin{theorem}\label{maintheorem}
Assume $(V_1)-(V_2)$, $(h_0)$ and $(f_0),(f_2)-(f_6)$. Then, there are $\epsilon_0, \mu^*$ and $t_1>0$, such that problem $(P)$ possesses a solution $v_\epsilon \in E $, with $I_\epsilon(v_\epsilon)=d_\epsilon >0$, for all $\epsilon \in (0,\epsilon_0),$ $t_0\in [0,t_{1})$ and $\mu \geq \mu^*$. Moreover, decreasing $\epsilon_0$ and $t_1$, and   increasing $\mu^*$, if necessary, we have two solutions $u_\epsilon, v_\epsilon \in E$ with
$$
I_\epsilon(u_\epsilon)=c_\epsilon<0<d_\epsilon=I_\epsilon(v_\epsilon).
$$ 
\end{theorem}

In the proof of Theorem \ref{maintheorem}, we use variational methods for nondifferentiable functional. A solution is obtained by applying Ekeland's variational principle, while the other one is obtained by using Mountain Pass Theorem.  Here, we would like point out that by applying the above theorem for the function $f$ given in (\ref{funcao}), we find two solutions $u_1,u_2 \in H^{1}(\mathbb{R}^{2})$ for the equation
$$
-\Delta{u}=2H(u-a)ue^{u^{2}}+\epsilon h, \quad \mbox{in} \quad \mathbb{R}^{2},
$$
with
$$
|[u_i=a]|=0 \quad \mbox{for} \quad i=1,2.
$$
\hspace{0.5 cm}

\noindent \textbf{Notation:} In this paper, we use the following
notations:
\begin{itemize}
	\item  The usual norms in $L^{t}(\mathbb{R}^{2})$ and $H^{1}(\mathbb{R}^{2})$ will be denoted by
	$|\,. \,|_{t}$ and $\|\;\;\;\|$ respectively.

	\item   $C$ denotes (possible different) any positive constant.
	
	\item   $B_{R}(z)$ denotes the open ball with center at $z$ and
	radius $R$.
	
	\item  If $B \subset \mathbb{R}^{2}$ is a mensurable set, let us denote by $|B|$ the Lebesgue's measure of $B$.
	
	\item $\Phi$ denotes the N-function $\Phi(t)=e^{|t|^{2}}-1$.
\end{itemize}

\section{Technical results involving the exponential critical growth}

 In this section, we will prove some technical lemmas, which are crucial in our approach. Since we will work with exponential critical growth, some versions of the Trudinger-Moser inequality are very important in our arguments. The first version that we would like to recall is due to Trudinger and Moser, see \cite{M} and \cite{T}, which claims if $\Omega$ is a bounded domain with smooth boundary, then for any  $ u \in H^{1}_0(\Omega)$,
\begin{equation} \label{X0}
\int_{\Omega}
e^{\alpha |u|^2}dx
< +\infty, \,\,\,\, \mbox{ for every }\,\,\alpha >0.
\end{equation}
Moreover, there exists a positive constant $C=C(\alpha,|\Omega|)$ such that
\begin{equation} \label{X1}
\sup_{||u||_{H_0^{1}(\Omega)} \leq 1} \int_{\Omega} e^{\alpha |u|^2} dx \leq C , \,\,\,\,\,\,\, \forall \, \alpha  \leq 4 \pi .
\end{equation}

A version in $H^{1}(\Omega)$ has been proved by Adimurthi and Yadava \cite{AY}, and it says that if $\Omega$ is a bounded domain with smooth boundary, then for any  $ u \in H^{1}(\Omega)$,
\begin{equation} \label{X02}
\int_{\Omega}
e^{\alpha |u|^2}dx
< +\infty, \,\,\,\, \mbox{ for every }\,\,\alpha >0.
\end{equation}
Furthermore, there exists a positive constant $C=C(\alpha,|\Omega|)$ such that
\begin{equation} \label{X11}
\sup_{||u||_{H^{1}(\Omega)} \leq 1} \int_{\Omega} e^{\alpha |u|^2} dx \leq C , \,\,\,\,\,\,\, \forall \, \alpha  \leq 2 \pi .
\end{equation}

The third version that we will use is due to Cao \cite{Cao}, which is version of the Trudinger-Moser inequality in whole space $\mathbb{R}^2$ and has the following statement:
\begin{equation} \label{X2}
\int_{\mathbb{R}^2}\left(e^{\alpha |u|^2}-1\right)dx<+\infty,\ \ \mbox{for all} \,\,\, u\in H^1(\mathbb{R}^2)\ \ \mbox{and} \,\, \alpha>0.
\end{equation}
Besides, given $\alpha < 4\pi$ and $M>0$, there is 
a constant $C_1=C_1(M,\alpha)>0$ verifying
\begin{equation} \label{X4}
\sup_{ u \in \mathcal{B}_M}
\int_{\mathbb{R}^2}\left(e^{\alpha |u|^2}-1\right)dx\leq C_1
\end{equation}
where
$$
\mathcal{B}_M=\{u \in H^{1}(\mathbb{R}^{2})\,:\, |\nabla u|_{2}\leq 1 \quad \mbox{and} \quad |u|_{2}\leq M  \}.
$$

As a consequence from (\ref{X2})-(\ref{X4}), we are able to prove some technical lemmas. The first of them is crucial in the study of the $(PS)$ condition for $I_\epsilon$. 
\begin{lemma} \label{alphat11} Let $\alpha>0$ and  $(u_{n})$ be a sequence in $H^{1}(\mathbb{R}^{2})$ with 
	$$
	\limsup_{n \to +\infty} \|u_n \|  < \sqrt{\frac{4 \pi}{\alpha}}.
	$$ 
	Then, there exist  $t> 1$, $t$ close to 1,  and $C > 0$  satisfying
	\[
	\int_{\mathbb{R}^{N}}\left(e^{\alpha |u_n|^{2}} - 1 \right)^t dx \leq C, \,\,\,\,\forall\, n \in \mathbb{N}. 
	\]

\end{lemma}
\noindent {\bf Proof.} \, 
As
\[
\limsup_{n\rightarrow\infty} \|u_n\| < \sqrt{\frac{4 \pi}{\alpha}},
\]
there are $m >0$ and $n_0\in \mathbb{N}$ verifying 
\[
\|u_n\|^{2} < m < \frac{4 \pi}{\alpha},
\,\,\,\,\forall\,  n \geq n_0.
\]
Fix  $t > 1$, with $t$ close to $1$, and $ \beta> t$ satisfying $\beta m < \frac{4 \pi}{\alpha}. $ Then, there exists $C=C(\beta)>0$ such that
$$
\int_{\mathbb{R}^{2}}\left(e^{\alpha |u_n|^{2}} -1 \right)^t dx   
\leq C\int_{\mathbb{R}^{2}}\left(e^{\alpha \beta m (\frac{|u_n|}{\|u_n\|})^{2}} - 1 \right) dx,
$$
for every  $n \geq n_0$.  Hence, by (\ref{X4}),
$$
\int_{\mathbb{R}^{2}}\left(e^{\alpha |u_n|^{2}} -1 \right)^t dx \leq C_1 \, \,\,\,\,\forall\,  n \geq n_0,
$$
for some positive constant $C_1$. Now, the lemma follows fixing
$$
C=\max\left\{C_1, \int_{\mathbb{R}^{2}}\left(e^{\alpha |u_{1}|^{2}} -1 \right)^t dx,....,\int_{\mathbb{R}^{2}}\left(e^{\alpha |u_{n_0}|^{2}} -1 \right)^t dx \right\}.
$$
\fim

\begin{lemma} \label{lemmaM}
Let $\beta,M>0$ verifying $\beta M < 4\pi$ and $q>2$. If $\parallel u\parallel^{2}\leq M$, then there is $C=C(\beta,M,q)>0$ such that,
$$
\int_{\mathbb{R}^2}|u|^q\left(e^{\beta |u|^2}-1\right)dx\leq C(\beta)\parallel u\parallel^q.
$$
\end{lemma}
\noindent {\bf Proof.} \, In what follows, fix $t>1$ close to 1, such that $\alpha=t \beta M < 4 \pi$. Then, there is  a constant $C>0$ such that
$$
\int_{\mathbb{R}^2}\left(e^{\beta |u|^2}-1\right)^{t}\,dx \leq C\int_{\mathbb{R}^2}\left(e^{\beta t |u|^2}-1\right)\,dx.
$$
Note that
$$
\int_{\mathbb{R}^2}\left(e^{\beta t |u|^2}-1\right)=\int_{\mathbb{R}^2}\left(e^{\beta t \|u\|^{2} (\frac{u}{\|u\|})^2}-1\right)\leq \int_{\mathbb{R}^2}\left(e^{\beta t M (\frac{u}{\|u\|})^2}-1\right)\,dx.
$$
Thereby, by (\ref{X4}),
$$
\int_{\mathbb{R}^2}\left(e^{\beta t |u|^2}-1\right) \leq \int_{\mathbb{R}^2}\left(e^{\alpha (\frac{u}{\|u\|})^2}-1\right)\,dx \leq \sup_{\|v\|\leq 1}\int_{\mathbb{R}^2}\left(e^{\alpha |v|^2}-1\right)\,dx=C<+\infty.
$$
From this, the function $\zeta_u=e^{\beta  |u|^2}-1 \in L^{t}(\mathbb{R}^{2})$ and there is $C=C(\beta,M)>0$ such that
\begin{equation} \label{limitacao}
|\zeta_u|_t \leq C, \quad \forall u \in B_M=\{u \in H^{1}(\mathbb{R}^{2})\,:\,\|u\|\leq M\}.
\end{equation}
Then, by applying the H\"older's inequality, 
$$
\int_{\mathbb{R}^2}|u|^q\left(e^{\beta |u|^2}-1\right)dx=\int_{\mathbb{R}^{2}}|u|^{q}\zeta_u\,dx \leq |\zeta_u|_t |u|_{t'q}^{q} 
$$ 
where $\frac{1}{t}+\frac{1}{t'}=1$. Hence, by Sobolev embedding, there is $C>0$ such that
\begin{equation} \label{limitacao2}
\int_{\mathbb{R}^2}|u|^q\left(e^{\beta |u|^2}-1\right)dx \leq C |\zeta_u|_t \|u\|^{q}. 
\end{equation}
Now, the lemma follows combining (\ref{limitacao}) and (\ref{limitacao2}). \fim

\section{Preliminaries about Orlicz spaces}

In this section, we recall some properties of Orlicz and Orlicz-Sobolev spaces. We refer to \cite{Adams, DT, FIN, Rao} for the fundamental properties of these spaces.
First of all, we recall that a continuous function $A : \mathbb{R} \rightarrow [0,+\infty)$ is a
N-function if:
\begin{description}
  \item[$(i)$] $A$ is convex.
  \item[$(ii)$] $A(t) = 0 \Leftrightarrow t = 0 $.
  \item[$(iii)$] $\displaystyle\lim_{t\rightarrow0}\frac{A(t)}{t}=0$ and $\displaystyle\lim_{t\rightarrow+\infty}\frac{A(t)}{t}= +\infty$ .
  \item[$(iv)$] $A$ is even.
\end{description}
We say that a N-function $A$ verifies the $\Delta_{2}$-condition, denote by  $A \in \Delta_{2}$, if
\[
A(2t) \leq K_*A(t)\quad \forall t\geq 0,
\]
for some constant $K_* > 0$. 

The complementary function ( or conjugate function ) $\widetilde{A}$ associated with $A$ is given
by the Legendre's transformation, that is,
\[
\widetilde{A}(s) = \max_{t\geq 0}\{ st - A(t)\} \quad  \mbox{for} \quad s\geq0.
\]
The functions $A$ and $\widetilde{A}$ are complementary each other. Moreover, we also have a Young type inequality given by
\begin{equation} \label{young}
st \leq A(t) + \widetilde{A}(s)\quad \forall s, t\geq0.
\end{equation}

In what follows, fixed an open set $\Omega \subset \mathbb{R}^{N}$ and a N-function $A$, we define the Orlicz space associated with $A$ as
\[
L^{A}(\Omega) = \left\{  u \in L_{loc}^{1}(\Omega) \colon \ \int_{\Omega} A\Big(\frac{|u|}{\lambda}\Big)dx < + \infty \ \ \mbox{for some}\ \ \lambda >0 \right\}.
\]
The space $L^{A}(\Omega)$ is a Banach space endowed with Luxemburg norm given by
\[
\Vert u \Vert_{A} = \inf\left\{  \lambda > 0 : \int_{\Omega}A\Big(\frac{|u|}{\lambda}\Big)dx \leq1\right\}.
\]
The convexity of $A$ implies in the inequality below, which will be used later on
\begin{equation} \label{P1}
\|u\|_A \leq 1 \Longleftrightarrow \int_{\Omega}A(|u|)dx \leq 1.
\end{equation}

Using the inequality (\ref{young}), it is possible to prove a H\"older type inequality, that is,
\[
\Big| \int_{\Omega}uvdx \Big| \leq 2 \Vert u \Vert_{A}\Vert v \Vert_{\widetilde{A}}\quad \forall u \in L^{A}(\Omega) \quad \mbox{and} \quad \forall v \in L^{\widetilde{A}}(\Omega).
\]
The space $L^{A}(\Omega)$ is separable and reflexive when $A$ and $\widetilde{A}$ satisfy the $\Delta_{2}$-condition. Moreover, the $\Delta_{2}$-condition implies that
\[
u_{n}\to u \ \mbox{in} \ \ L^{A}(\Omega)\quad \Leftrightarrow \quad \int_{\Omega}A(\vert u_{n} - u\vert)dx \rightarrow 0.
\]

\subsection{ The class $K_A(\Omega)$ and the subspace $E_A(\Omega)$}

In the study of the Orlicz space $L^{A}(\Omega)$, we denote by $K_A(\Omega)$ the following set
$$
K_A(\Omega)=\left\{u:\Omega \to \mathbb{R}\,:\, u \, \mbox{is mensurable and} \, \int_{\Omega}A(|u|)\,dx\leq 1 \right\}
$$
and let us denote by $E_A(\Omega) \subset L^{A}(\Omega)$ the following subspace
$$
E_A(\Omega)=\overline{L^{\infty}(\Omega)}^{\|\,\,\|_A} \quad \mbox{if} \quad |\Omega|<+\infty \quad \mbox{is bounded}
$$
or
$$
E_A(\Omega)=\overline{C^{\infty}_0(\Omega)}^{\|\,\,\|_A} \quad \mbox{if} \quad |\Omega|=+\infty \quad \mbox{is unbounded}.
$$
Using the above notations, it follows that
$$
E_A(\Omega) \subset K_A(\Omega) \subset L^{A}(\Omega),  
$$
and
\begin{equation} \label{KA}
K_A(\Omega)\subset \{u \in L^{A}(\Omega)\,:\,dist(u,E_A(\Omega))\leq 1\}.
\end{equation}

It is possible to prove that if $A$ verifies $\Delta_2$-condition, then
$$
E_A(\Omega) = K_A(\Omega) = L^{A}(\Omega).
$$ 

However, if $A$ does not satisfies the  $\Delta_2$-condition, we have that $E_A(\Omega)$ is a proper subspace of $L^A(\Omega)$. For example, this situation holds for the N-function $\Phi(t)=e^{|t|^{2}}-1$, because it does not verify the $\Delta_2$-condition. Moreover,  $L^{\Phi}(\Omega)$ is not reflexive, hence we cannot guarantee that if $J_0 \in (L^{\Phi}(\Omega))^*$, then
$$
J_0(u)=\int_{\mathbb{R}^{2}}vu\,dx, \forall u \in L^{\Phi}(\mathbb{R}^2),
$$ 
for some mensurable function $v:\mathbb{R}^{2} \to \mathbb{R}$. However, this type of problem does not hold in $(E_{\Phi}(\Omega))^*$, because if $J_1 \in (E_{\Phi}(\Omega))^*$ we know that there exists $v \in L^{\widetilde{\Phi}}(\mathbb{R}^2)$ such that
$$
J_1(u)=\int_{\mathbb{R}^{2}}vu\,dx, \forall u \in L^{\Phi}(\mathbb{R}^2).
$$

\begin{lemma}\label{LtildePhi}
Let $\xi(t)=\max\{t,t^2\}$ and $\widetilde{\Phi}$ the conjugate function associated with $\Phi$. Then, 
$$
\tilde\Phi\left(\frac{\Phi(r)}{r}\right )\leq \Phi(r) \quad \text{ and } \quad \tilde \Phi(t r)\leq \xi(t)\tilde \Phi (r), \ t,r\geq 0.
$$
Hence, $\tilde{\Phi}\in \Delta_2$, $E_{\tilde \Phi}(\mathbb{R}^2)=L_{\tilde\Phi}(\mathbb{R}^2)$ and $L_{\tilde\Phi}(\mathbb{R}^2)$ is separable. 
\end{lemma}
{\bf Proof.} \, The first inequality follows from \cite{FIN}. To prove the second one, we recall that  
$$
2\leq \frac{\Phi'(t)t}{\Phi(t)}, \ t\in (0,+\infty).
$$
Fix  $s>0$, such that $t=\tilde \Phi'(s)$. Since $\tilde\Phi'=(\Phi')^{-1}$ and $s\tilde\Phi'(s)=\tilde \Phi(s)+\Phi(\tilde \Phi'(s))$, we derive that 
$$
2\leq \frac{\Phi'(\tilde\Phi'(s))\tilde \Phi'(s)}{\Phi(\tilde \Phi'(s))}=\frac{s\tilde \Phi'(s)}{s\tilde\Phi'(s)-\tilde \Phi(s)},
$$
and so, 
$$
2s\tilde\Phi'(s)-2\tilde \Phi(s)\leq s\tilde \Phi'(s),
$$
that is
$$
\frac{s\tilde\Phi'(s)}{\tilde\Phi(s)}\leq 2.
$$
Now, fixing $s=\rho r>0$, we get
$$
\frac{d}{d\sigma}\left(\ln \left(\tilde\Phi(\rho r)\right)\right)\leq \frac{2}{\sigma}.
$$
From this, 
\begin{equation}\label{Pri1des}
\tilde \Phi(t r)\leq t^2\tilde \Phi (r), \ t\geq 1 \text{ and }r\geq 0.
\end{equation}
On the other hand, the convexity of $\tilde\Phi$ combines with $\tilde \Phi(0)=0$ to give 
$$
1\leq \frac{\tilde \Phi'(t)t}{\tilde\Phi(t)}, \ t\in (0,+\infty).
$$
Using again \cite{FIN}, we see that 
\begin{equation}\label{Seg2des}
\tilde \Phi(t r)\leq t\tilde \Phi (r), \ t\in(0,1] \text{ and }r\geq 0.
\end{equation}
Hence, from (\ref{Pri1des}) and (\ref{Seg2des}),
$$
\tilde \Phi(t r)\leq \xi(t)\tilde \Phi (r), \ t,r\geq 0. 
$$
Now, the conclusion follows from  \cite{FIN}.  $\hfill{\rule{2mm}{2mm}}$

\begin{lemma}\label{imeEPHI}
Let $X=H_0^1(\Omega)$ or $X=H^1(\Omega)$, where $\Omega \subset \mathbb{R}^2$ is a smooth bounded domain or $\Omega=\mathbb{R}^{2}$. Then, the embedding  $X\hookrightarrow E_\Phi(\Omega)$ is continuous. 
\end{lemma}
{\bf Proof.} \, From (\ref{X1}),(\ref{X11}) and (\ref{X4}), we have that the embedding $ X\hookrightarrow L_\Phi(\Omega)$ is continuous. Then, the lemma follows by demonstrating the inclusion $X\subset E_\Phi(\Omega)$. First of all, by  (\ref{X0}),(\ref{X02}) and (\ref{X2}), we know that    
$$
\int_\Omega (e^{\lambda |u|^2}-1)dx<+\infty, \quad \forall \lambda \geq 0 \quad \mbox{and} \quad  \forall u\in X, 
$$
implying that 
\begin{equation}\label{l1}
X\subset K_\Phi(\Omega).
\end{equation}
Assume by contradiction that there is $u_0\in X$ with $u_0\not\in E_\Phi(\Omega)$. Since $E_\Phi(\Omega)$  is a closed subspace in $L^\Phi(\Omega)$, we ensure that $dist(u_0,E_\Phi(\Omega))>0$. Thereby, for 
$$
\lambda>\frac{1}{dist(u_0,E_\Phi(\Omega))},
$$
we have that
\begin{eqnarray*}
dist(\lambda u_0,E_\Phi(\Omega))&=&\inf\left\{\|\lambda u_0-v\|_\Phi; v\in E_\Phi(\Omega)\right\}\\
&=&\lambda\inf\left\{\| u_0-\frac{v}{\lambda}\|_\Phi; \frac{v}{\lambda}\in E_\Phi(\Omega)\right\}\\
&=&\lambda dist(u_0,E_\Phi(\Omega))>1.
\end{eqnarray*}
Then, by (\ref{KA}), $\lambda u_0\not\in K_\Phi(\Omega)$, which contradicts (\ref{l1}), because $\lambda u_0 \in X$.
$\hfill{\rule{2mm}{2mm}}$

\begin{lemma}\label{imer}
The embeddings $E_\Phi(\Omega)\hookrightarrow L^{2n}(\Omega)$ are continuous for any $n\in\mathbb{N}$.
\end{lemma}
{\bf Proof.} \, For each $n\in \mathbb{N}^*$, we know that 
$$
\frac{1}{n!}t^{2n}\leq \sum_{k=0}^{+\infty}\frac{1}{k!}t^{2k}=e^{t^2}-1,\quad \forall t\in\mathbb{R}.
$$
Then, for each $u\in E_\Phi(\Omega)$
$$
\frac{1}{n!}\int_\Omega \left(\frac{u}{|u|_\Phi}\right)^{2n}dx\leq \int_\Omega\left(e^{\left(\frac{u}{|u|_\Phi}\right)^2}-1\right)dx\leq 1,
$$
leading to 
$$
|u|^{2n}_{2n}=\int_\Omega |u|^{2n}\leq n!|u|_\Phi^{2n},
$$
showing the lemma. $\hfill{\rule{2mm}{2mm}}$

\section{Some properties of the functional $\Psi$}

Let $f:\mathbb{R}^2\times \mathbb{R}\to \mathbb{R}$  be a mensurable function for each  $t\in \mathbb{R}$ and Locally Lipschitzian  for each $x\in\mathbb{R}^2$ verifying:

\begin{description}
  \item[($f_2$)] There is $t_0\geq0$ such that 
$$
f(x,t)=0  \quad \text{for } \quad  t< t_0  \text{ and } \, \forall x \in \mathbb{R}^{2}
$$
and
$$
f(x,t)>0 \quad \mbox{for} \quad  t>t_0 \text{ and } \, \forall x \in \mathbb{R}^{2}.
$$
\item[($f_*$)] There are $\alpha_0,c_1,c_2>0$ and $\alpha_0>0$ such that 
$$
\left|\xi\right|\leq c_1|u|+c_2e^{\alpha_0|u|^2}, \ \forall \xi\in \partial_t F(x,u) \quad \mbox{and} \quad \forall x\in\mathbb{R}^2,
$$
where $F(x,t)=\int_0^t f(x,s)ds$.
\end{description}

\begin{theorem}\label{Liploc}
The functional $\Psi:E_\Phi(\Omega)\to \mathbb{R}$ given by 
$$
\Psi(u)=\int_\Omega F(x,u)dx
$$
is well defined and $\Psi\in Lip_{loc}(E_{\Phi}(\Omega),\mathbb{R})$.
\end{theorem}
{\bf Proof.} \, For each $u\in E_\Phi(\Omega)$ and $R>0$, consider $w,v\in B_R(u)\subset E_\Phi(\Omega)$. By Lebourg's Theorem, there is $\xi\in\partial F_t(x,\theta)$ with $\theta\in [w,v]$ such that 
$$
\left|F(x,w)-F(x,v)\right|=|\langle\xi,w-v\rangle|\leq |\xi||w-v|.
$$
Then by  $(f_*)$, 
$$
\left|F(x,w)-F(x,v)\right|\leq\left ( c_1|\theta|+c_2(e^{\alpha|\theta|^2}-1)\right)|w-v|,
$$
for $\alpha > \alpha_0$ and $\alpha$ close to $\alpha_0$. Setting $\eta(x)=|v(x)|+|w(x)|$, it follows that 
$$
\left|F(x,w)-F(x,v)\right|\leq  \left(c_1|\eta|+c_2(e^{\alpha|\eta|^2}-1)\right)|w-v|,
$$
and so,
$$
\left|\Psi(w)-\Psi(v)\right|\leq \int_\Omega \left(c_1|\eta|+c_2(e^{\alpha|\eta|^2}-1)\right)|w-v|\,dx.
$$
By H\"older's inequality,
$$
\left|\Psi(w)-\Psi(v)\right|\leq c_1|\eta|_2|w-v|_2+c_2\left(\int_\Omega (e^{2\alpha|\eta|^2}-1)dx\right)^{\frac{1}{2}}|w-v|_2.
$$
Once $E_\Phi(\Omega)\hookrightarrow L^{2}(\Omega)$ is continuous, see Lemma \ref{imer}, we derive that
\begin{eqnarray}\label{d0}
\left|\Psi(w)-\Psi(v)\right|&\leq& c_1\left( |w-u|_\Phi+|u-v|_\Phi+2|u|_\Phi\right)|w-v|_\Phi\nonumber\\
&&+c_2\left(\int_\Omega (e^{2\alpha|\eta|^2}-1)dx\right)^{\frac{1}{2}}|w-v|_\Phi .
\end{eqnarray}
On the other hand, the convexity of $\Phi$ yields 
\begin{eqnarray}\label{d1}
\int_\Omega \left( e^{2\alpha(|w|+|v|)^2}-1\right)dx&\leq&
\frac{1}{4}\int_\Omega \left( e^{32\alpha|w-u|^2}-1\right)dx
+\frac{1}{4}\int_\Omega \left( e^{32\alpha|v-u|^2}-1\right)dx\nonumber\\
&&+\frac{1}{2}\int_\Omega \left( e^{32\alpha|u|^2}-1\right)dx.
\end{eqnarray}
Now, fixing $R>0$ verifying $R <\frac{1}{\sqrt{32\alpha_0}}$ and $\alpha$ close to $\alpha_0$ satisfying $R<\frac{1}{\sqrt{32\alpha}}$, we derive that
$$
|\sqrt{32\alpha}(w-u)|_\Phi\leq \sqrt{32\alpha}R\leq 1 \text{ and }|\sqrt{32\alpha}(v-u)|_\Phi\leq \sqrt{32\alpha}R\leq 1,
$$
and so,
\begin{equation}\label{d2}
\int_\Omega \left( e^{32\alpha|w-u|^2}-1\right)dx\leq1 \quad \text{and } \quad \int_\Omega \left( e^{32\alpha|v-u|^2}-1\right)dx\leq1,
\end{equation}
for all  $w,v\in B_{R}(u)$. From (\ref{d1})-(\ref{d2})
\begin{equation}\label{d3}
\int_\Omega \left( e^{2\alpha(|w|+|v|)^2}-1\right)dx\leq \frac{1}{4}\left(2+2\int_\Omega \left( e^{32\alpha|u|^2}-1\right)dx\right).
\end{equation}
Thereby, gathering (\ref{d0}) and (\ref{d3}),
\begin{eqnarray*}
\left|\Psi(w)-\Psi(v)\right|&\leq& c_1\left( R+R+2|u|_\Phi\right)|w-v|_\Phi\nonumber\\
&&+c_2 \frac{1}{4}\left(2+2\int_\Omega \left( e^{32\alpha|u|^2}-1\right)dx\right)^{\frac{1}{2}}|w-v|_\Phi\\
&=&K(R,u)|w-v|_\Phi, \ \forall w,v\in B_{R}(u).
\end{eqnarray*}
$\hfill{\rule{2mm}{2mm}}$

Our next goal is proving the differentiable inclusion  
\begin{equation} \label{Inc}
\partial \Psi(u)\subset \int_\Omega \partial_t F(x,u)dx, \ u\in E_\Phi(\Omega).
\end{equation}
To do this, we need of the following result
\begin{lemma}\label{convE}
Let $\psi:\mathbb{R}\to\mathbb{R}_+$ be a N-function and 
$$
g_n\to g \text{ in }E_\psi(\Omega).
$$
Then, there is $\hat{g} \in E_\psi(\Omega)$ and a subsequence of $\{g_n\}$, denoted by $\{g_{m_k}\}$, such that 
\begin{description}
\item{$(i)$} $g_{m_k}(x)\to g(x)\quad  \text{ a.e. } \quad   x\in \Omega,$
\item{$(ii)$} $|g_{m_k}(x)|\leq \hat{g}(x)\quad \text{ a.e. } \quad x\in \Omega.$
\end{description}
\end{lemma}
{\bf Proof.} \, As 
$$
|g_m-g|_\psi\to 0,
$$
we have that
$$
\int_\Omega \psi(g_m-g)dx\leq |g_m-g|_\psi \to 0,
$$
implying that there is a subsequence of $\{g_n\}$, denoted by $\{g_{m_k}\}$, such that 
$$
\psi(g_{m_k}-g)(x)\to 0 \quad \text{ a.e. in } \quad \Omega,
$$
and so, 
$$
(g_{m_k}-g)(x)=\psi^{-1}\circ\psi(|g_{m_k}-g|)(x)\to 0 \quad \text{ a.e. in } \quad \Omega,
$$
that is,
$$
g_{m_k}(x)\to g(x) \quad \text{a.e. in} \quad \Omega,
$$
Now, define 
$$
\zeta_m=\displaystyle \sum_{k=1}^{m}|g_{n_{k+1}}-g_{n_k}|\in E_\psi(\Omega),
$$
with
$$
|g_{n_{k+1}}-g_{n_k}|_\psi<\frac{1}{2^k}, \quad  \forall k\in\mathbb{N}.
$$
Hereafter, $g_k$ denotes $g_{n_k}$, that is, $g_k:=g_{n_k}$. For $n\leq m$, 
$$
|\zeta_m-\zeta_n|_{\psi}\leq \sum_{k=n}^{m}|g_{{k+1}}-g_{k}|_{\psi} \leq \sum_{k=n}^{m}\frac{1}{2^k}\to 0, 
$$
from if follows that $\{\xi_m\}\subset E_\psi(\Omega)$ is a Cauchy's sequence in $E_\psi(\Omega)$. Once $E_\psi(\Omega)$ is a Bancah space, there exists $\zeta \in E_\psi(\Omega)$ such that 
$$
\zeta_m\to \zeta \text{ in } E_\psi(\Omega).
$$
Then
$$
\int_\Omega\psi(\xi_m-\xi)dx\leq |\zeta_m-\zeta|_\psi\to 0,
$$
and so, 
$$
\zeta_{m_k}(x)\to \zeta(x) \text{ a.e. in } \Omega,
$$
and
$$
\zeta_{m_k}(x)\leq \zeta(x) \text{ a.e. in } \Omega, \ k\in\mathbb{N}.
$$
On the other hand, for $n\leq m$,
$$
|g_m-g_n|(x)\leq \xi_m(x)\leq \zeta(x)\text{ a.e. in } \Omega.
$$
Setting $\hat{g}=\zeta+|g| \in E_\psi(\Omega)$ and taking the  $n\to+\infty$, we get 
$$
|g_m(x)|\leq \hat{g}(x)  \text{ a.e. in } \Omega \quad \forall  m\in\mathbb{N},
$$
showing  $(ii)$. $\hfill{\rule{2mm}{2mm}}$
\begin{theorem}\label{theoInclusao}
Assume $(f_*)$ and that $\underline{f}(x,t)$ and $\overline{f}(x,t)$ are N- mensurable functions. If $\Omega \subset \mathbb{R}^2$ is a smooth bounded domain or $\Omega=\mathbb{R}^{2}$, then for each  $u\in E_{\Phi}(\Omega)$,  
\begin{equation}\label{inclusao}
\partial \Psi(u)\subset \partial_t F(x,u)=[\underline{f}(x,u(x)),\overline{f}(x,u(x))] \,  \text{a.e. in }\Omega.
\end{equation}
Moreover,
$$
\partial \Psi|_X(u)\subset \partial \Psi(u), \ u\in X,
$$
where $X=H_0^1(\Omega)$ or $X=H^1(\Omega)$. Here, the above inclusion means that given $\xi\in \partial \Psi(u)\subset E_\Phi(\Omega)^*$, there is  $\tilde{\xi}\in L^{\tilde{\Phi}}(\Omega)$ satisfying 
\begin{itemize}
\item $\langle\xi,v\rangle=\int_\Omega \tilde{\xi}vdx, \quad \forall v\in E_\Phi(\Omega)$,
\item $ \tilde{\xi}(x)\in \partial_tF(x,u)=[\underline{f}(x,u(x)),\overline{f}(x,u(x))] \text{ a.e. in }\Omega.$
\end{itemize}
\end{theorem}
{\bf Proof.} \, Given $u,v\in E_\Phi(\Omega)$, let $\{g_j\}\subset E_\Phi(\Omega)$ with $ g_j\to 0$ in $E_\Phi(\Omega)$ and
$\{\lambda_j\}\subset \mathbb{R}_+$ with $\lambda_j\to 0$ verifying
\begin{equation}\label{limsup}
\Psi^0(u;v)=\displaystyle\lim_{j\to+\infty}\int_\Omega \frac{F(u+g_j+\lambda_jv)-F(u+g_j)}{\lambda_j}dx.
\end{equation}
Setting
$$
F_j(u;v):=\frac{F(u+g_j+\lambda_jv)-F(u+g_j)}{\lambda_j},
$$
the Lebourg's Theorem guarantees that there is $\xi_j\in \partial_tF(x,\theta_j)$, with $\theta_j\in[u+g_j+\lambda_jv,u+g_j]$ such that 
$$
\left|F_j(u;v)\right|=\frac{1}{\lambda_j}|\langle \xi_j,\lambda_j v\rangle|\leq |\xi_j||v|.
$$
Hence by $(f_*)$, 
$$
\left|F_j(u;v)\right|\leq \left(c_1|\theta_j|+c_2(e^{\alpha|\theta_j|^2}-1)\right)|v|,
$$
for $\alpha > \alpha_0$ and $\alpha$ close to $\alpha_0$. Fixing 
$$
\beta_j=(|u|+|g_j|+\lambda_j|v|)+(|u|+|g_j|)=2|u|+2|g_j|+\lambda_j|v|,
$$
we see that
\begin{equation}\label{desFj}
\left|F_j(u;v)\right|\leq \left(c_1|\beta_j|+c_2(e^{\alpha|\beta_j|^2}-1)\right)|v|.
\end{equation}
Applying Lemma \ref{convE}, there exists $g_*\in E_\Phi(\Omega)$ such that 
\begin{equation}\label{desbeta}
|\beta_j|\leq 2|u|+2g_*+c|v|\, \text{a.e. in }\Omega,
\end{equation}
for some subsequence. Thereby, by (\ref{desFj}) and (\ref{desbeta}), there exists a subsequence $\{F_{j_k}(u;v)\}$ such that 
$$
\left|F_{j_k}(u;v)\right|\leq \left(c_1(2|u|+2g_*+c|v|)+c_2(e^{\alpha(2|u|+2g_*+c|v|)^2}-1)\right)|v|\in L^1(\Omega).
$$
Applying the Lebesgue's Theorem, 
\begin{eqnarray}
\Psi^0(u;v)&=&\displaystyle \lim_{j_k\to+\infty} \int_\Omega F_{j_k}(u;v)dx= \int_\Omega\lim_{j_k\to+\infty} F_{j_k}(u;v)dx\nonumber\\
&\leq& \int_\Omega F^0(u;v)dx=\int_\Omega \max\{\langle\xi,v\rangle; \xi\in\partial_t F(x,u)\}dx\nonumber\\
&\leq& \int_{[v<0]}\underline{f}(x,u)vdx+ \int_{[v>0]}\overline{f}(x,u)vdx.
\end{eqnarray}
Now, we will show that for each $\xi\in\partial \Psi(u)\subset (E_\Phi(\Omega))^*$, the function $\tilde{\xi}\in L^{\widetilde{\Phi}}(\Omega)$, which satisfies  
$$
\langle \xi,w\rangle=\int_\Omega \tilde{\xi}wdx, \quad  \forall w\in E_\Phi(\Omega),
$$
must verify
$$
\tilde{\xi}(x)\in [\underline{f}(x,u(x)),\overline{f}(x,u(x))] \quad \text{a.e. in} \quad \Omega.
$$
Indeed, assume by contradiction that there is a mensurable set $\mathcal{M} \subset \Omega$, with $0<|\mathcal{M}|<+\infty$, satisfying 
\begin{equation}\label{desf1}
\tilde{\xi}(x)<\underline{f}(x,u(x)), \ x\in \mathcal{M}.
\end{equation}
Setting $v=-\chi_{\mathcal{M}} \in E_\Phi(\Omega)$, we must have  
$$
-\int_{\mathcal{M}}\tilde{\xi}dx=\int_\Omega\tilde{\xi}\left(-\chi_\mathcal{M}\right)dx\leq \Psi^0(u,-\chi_{\mathcal{M}})\leq-\int_{\mathcal{M}} \underline{f}(x,u(x))dx,
$$
leading  to
$$
\int_\Omega\tilde{\xi}\chi_{\mathcal{M}}dx\geq \int_{\mathcal{M}} \underline{f}(x,u(x))dx,
$$
which contradicts (\ref{desf1}). Thereby,
$$
\tilde{\xi}(x) \geq \underline{f}(x,u(x)) \quad \mbox{a.e in} \quad \Omega.
$$
The same type of arguments work to show that  
$$
\tilde{\xi}(x)\leq\overline{f}(x,u(x))\, \text{a.e. in } \Omega.
$$

From definition of $X$, we know that  $\overline{X}^{\|\,\,\,\|_\Phi}=E_\Phi(\Omega) $, then the Lemma \ref{imeEPHI} combined with chain rule gives 
$$
\partial \Psi|_X(u)\subset \partial \Psi(u), \quad  \forall u\in X.
$$
$\hfill{\rule{2mm}{2mm}}$

\section{An aplication}
In this section, we will study the existence of solution for the following class of multivalued elliptic equation
$$
-\Delta u+V(x)u-\epsilon h(x)\in \partial_t F(x,u), \quad \text{in} \quad \mathbb{R}^2, \eqno{(P)}
$$
where
\begin{itemize}
\item $h\in H^{-1}$, that is, the functional $\langle h,v\rangle=\displaystyle \int_{\mathbb{R}^2} hvdx$ is continuous in $H^1(\mathbb{R}^2) $ and $0< \displaystyle \int_{\mathbb{R}^{2}}h\,dx< +\infty$. 
\item $F(x,t) =\displaystyle \int_0^tf(x,s)ds$, for $(x,t)\in\mathbb{R}^2\times\mathbb{R}$, where $f$ verifies $(f_0)$ and $(f_1)$.
\end{itemize}
Related to the potential $V:\mathbb{R}^2\to \mathbb{R}$, we assume that 
\begin{description}
\item[$(V_1)$] $V$ is continuous and $V(x)\geq V_0>0$, $\forall x\in\mathbb{R}^2$,
\item[$(V_2)$] $\frac{1}{V}\in L^1(\mathbb{R}^2)$.
\end{description}
In order to apply variational methods, we will consider the Hilbert space 
$$
E:=\left\{u\in H^1(\mathbb{R}^2)/ \int_{\mathbb{R}^2}V(x)|u|^2dx<+\infty\right\}
$$
endowed with the inner product 
$$
\left\langle u,v \right\rangle=\int_{\mathbb{R}^{2}}(\nabla u \nabla v +V(x)uv)\,dx .
$$
Associated with the above inner product, we have the norm
$$
\parallel u\parallel=\left(\int_{\mathbb{R}^2} (|\nabla u|^2+V(x)|u|^2)\,dx\right)^{\frac{1}{2}} .
$$
using the above information, it is well known that
\begin{description}
\item [$(E_1)$] $E\hookrightarrow L^q(\mathbb{R}^2)$ is a compact embedding for all $q\geq 1$, see \cite{JEU1,JEU2}
\item [$(E_2)$] $E\hookrightarrow H^1(\mathbb{R}^2)\hookrightarrow E_\Phi(\Omega)$ is a continuous embedding (see Lemma \ref{imeEPHI}).
\end{description}
In the present paper, we say that $u\in E$ is a solution for $(P)$, if there is $\rho\in L_{\tilde\Phi}(\mathbb{R}^2)$ such that 
\begin{description}
\item[$(i)$] $\displaystyle \int_{\mathbb{R}^2}(\nabla u\nabla v+V(x)uv)dx-\int_{\mathbb{R}^2}\rho vdx-\epsilon\int_{\mathbb{R}^2} h vdx=0, \ v\in E$,
\item[$(ii)$] $\rho(x)\in \partial_t F(x,u(x))$ \quad \mbox{a.e. in} \quad $\mathbb{R}^2$,
 \item[$(iii)$] $\left|[u>t_0]\right|>0$.
\end{description}

The reader is invited to observe that  $u \in E$ is a solution for $(P)$ if, and only if, $u$ is a critical point of the energy functional associated with $(P)$ given by:
$$
I_\epsilon(u)=\int_{\mathbb{R}^2}(|\nabla u|^2+V(x)|u|^2)dx -\int_{\mathbb{R}^2} F(x,u)dx -\epsilon\int_{\mathbb{R}^2} hvdx, \ u\in E.
$$
Note that Theorem \ref{Liploc} gives that $I_\epsilon\in Lip_{loc}(E;\mathbb{R})$. From this, using some properties of the generalizity gradient together with Theorem \ref{theoInclusao}, given $u\in E$ and $w\in\partial I_\epsilon (u)$, there exists $\rho\in L_{\tilde\Phi}(\mathbb{R}^2)$ such that
$$
\langle w,v\rangle= \int_{\mathbb{R}^2}(\nabla u\nabla v + V(x)uv)dx-\int_{\mathbb{R}^2}\rho vdx-\epsilon \int_{\mathbb{R}^2} h vdx \quad \forall v\in E,
$$
with
$$
\rho(x)\in \partial_t F(x,u(x)) \quad \text{a.e. in } \mathbb{R}^2.
$$

\section{Existence of solution via Ekeland's variational principle}
In this section, we will get a solution via Ekeland's variational principle. 

\begin{lemma}\label{Inf} Assume that $(f_0)$ and $(f_2)-(f_4)$ hold. Then, there are $\epsilon_0,r, \alpha, \delta >0$, such that  
$$
c_\epsilon=\inf_{\|u\| \leq r}I_\epsilon(u) <-\delta
$$
and
$$
I_\epsilon(u) \geq \alpha \quad \mbox{for} \quad \|u\|=r
$$
for all $\epsilon \in (0, \epsilon_0)$. Here, $r$ is independent of $\epsilon$, but $\alpha$ and $\delta$ depend on $\epsilon$. Moreover, the numbers $\epsilon_0,r, \alpha$ and $\delta$ do not depend  on $t_0$ given in $(f_2)$. 
\end{lemma}
{\bf Proof.} Using the conditions on $F$, given $\beta \in (0, \lambda_1)$, $q>2$ and $\alpha > \alpha_0$ close to $\alpha_0$, we have that
$$
F(x,t)\leq \frac{(\lambda_1 -\beta)}{2}|t|^{2} +C|t|^{q}(e^{\alpha |t|^{2}}-1), \quad \forall t \in \mathbb{R}.
$$
Then, fixing $r>0$ small enough such that $\alpha r^{2} < 4 \pi$  and using Lemma \ref{lemmaM}, we get for $u \in E$ with $\|u\|\leq r$, 
\begin{eqnarray*}
I_\epsilon(u)&\geq& \frac{1}{2}\parallel u\parallel^2-\frac{(\lambda_1- \beta)}{2}|u|_2^2-C\parallel u\parallel^q-\epsilon\parallel h\parallel_*\parallel u\parallel\\
&=&\frac{1}{2}\Big(1-\frac{(\lambda_1-\beta)}{\lambda_1}\Big)\parallel u\parallel^2-C\parallel u\parallel^q-\epsilon \parallel h\parallel_*\parallel u\parallel,\\
\end{eqnarray*}
showing that $I_\epsilon$ is bounded from below  for $\|u\|\leq r$. Moreover, decreasing $r$ if necessary of a way that 
$$
\frac{1}{2}r^{2}-Cr^{q} \geq \frac{1}{4}r^{2},
$$
we derive that
$$
I_\epsilon(u) \geq \frac{1}{4}r^{2} - \epsilon \|h\|r, \quad \|u\|=r.
$$
Thereby, choosing $\epsilon_0>0$ such that
$$
\alpha_\epsilon=\frac{1}{4}r^{2} - \epsilon \|h\|r>0, \quad \forall \epsilon \in (0, \epsilon_0),
$$
we see that 
$$
I_\epsilon(u) \geq \alpha_\epsilon \quad \mbox{for} \quad \|u\|=r, \quad \forall \epsilon \in (0, \epsilon_0).
$$
Now, take $v \in E$ satisfying 
$$
\|v\|=1 \quad \mbox{and} \quad \int_{\mathbb{R}^{2}}hv\,dx>0.
$$
Note that for each $s>0$, 
$$
I_\epsilon(sv)=\frac{s^2}{2}-\int_{\mathbb{R}^2}F(x,sv)dx -\epsilon s\int_{\mathbb{R}^2} hvdx<\frac{s^2}{2}-\epsilon s\int_{\mathbb{R}^2} hvdx.
$$
Fixing $s=s(\epsilon)>0$ small enough  satisfying 
$$
\delta= -\frac{s^2}{2}+\epsilon s\int_{\mathbb{R}^2} hvdx>0,
$$
it follows that $\|sv\|<r $ and 
$$
I_\epsilon(sv)<-\delta <0,
$$
implying that 
$$
c_\epsilon=\displaystyle\inf_{\| u\|\leq r}I_\epsilon(u)<-\delta<0.
$$
$\hfill{\rule{2mm}{2mm}}$
\begin{theorem}\label{SolEkland}
Assume $(V_1),(V_2), (f_0),$ $(f_2)$ and $(f_3)$. Then, problem $(P)$ possesses a solution $u_\epsilon \in E $, with $I_\epsilon(u_\epsilon)=c_\epsilon<-\delta<0$, for all $\epsilon \in (0,\epsilon_0)$ and $t_0\in [0,t_{\ast})$, with $t_{\ast}=t_{\ast}(\epsilon)=\frac{2\delta}{\epsilon \int_{\mathbb{R}^2}hdx}>0$.
\end{theorem}
{\bf Proof.}\, Fix $r>0$ such that $\alpha_0 r^{2} < 4 \pi$. Applying the Lemma \ref{Inf} together with Ekeland's variational principle, there is $\{u_n\}\subset \overline{B}_{r}(0)$ verifying 
\begin{itemize}
\item $I_\epsilon(u_n)\to c_\epsilon$ (as $n\to +\infty$),
\item $\lambda_\epsilon(u_n):=\min\{\parallel \xi\parallel_{E^*}/\xi\in \partial I_\epsilon(u_n)\}\to 0$ (as $n\to +\infty$).
\end{itemize}
Next, we  fix $w_n \in \partial I_\epsilon (u_n)$ and $\{\rho_n\}\subset L_{\tilde\Phi}(\mathbb{R}^2)$ verifying  
$$
\parallel w_n\parallel_{E^*}:=\lambda_\epsilon(u_n)
$$
\begin{equation}\label{IG}
\langle w_n,v\rangle=\int_{\mathbb{R}^2}\nabla u_n\nabla v +V(x)u_n v dx-\int_{\mathbb{R}^2}\rho_n vdx-\epsilon \int_{\mathbb{R}^2}hvdx, \,\, \forall v\in E,
\end{equation}
and
$$
\rho_n(x)\in \partial_t F(x,u_n(x)) \quad \text{ a.e. in }  \, \mathbb{R}^2.
$$
We claim that $\{\rho_n\}$ is bounded in $L_{\tilde\Phi}(\mathbb{R}^2)$. Indeed, fixing $p>4$, $\alpha >\alpha_0$ with $\alpha r^{2}<4 \pi$, and using $(f_*)$, we get 
$$
\int_{\mathbb{R}^2}\tilde \Phi(\rho_n)dx\leq \int_{\mathbb{R}^2}\tilde \Phi\left(c_1|u_n|+c_2|u_n|^p(e^{\alpha|u_n|^2}-1)\right)dx.
$$
The convexity of $\tilde \Phi$ and the $\Delta_2$- condition combine to give 
$$
\int_{\mathbb{R}^2}\tilde \Phi(\rho_n)dx\leq \frac{\xi(2c_1)}{2}\int_{\mathbb{R}^2}\tilde \Phi(|u_n|)dx+\frac{\xi(2c_2)}{2}\int_{\mathbb{R}^2}\tilde\Phi\left(|u_n|^p(e^{\alpha|u_n|^2}-1)\right)dx.
$$
By Lemma \ref{LtildePhi} and $(E_1)$, there are positive constants $C_1,C_2$ such that
\begin{eqnarray*}
\int_{\mathbb{R}^2}\tilde \Phi(\rho_n)dx&\leq&C_1(|u_n|_1+|u_n|_2^2) + \\
&&+ C_2\int_{\mathbb{R}^2}(|u_n|^{p+1}+|u_n|^{2(p+1)})\left(e^{\alpha_0 |u_n|^2}-1\right)dx.
\end{eqnarray*}
Recalling that the space $E$ is continuously embedding in $L^1(\mathbb{R}^2)$ and $L^2(\mathbb{R}^2)$, and $\alpha r^{2} <4\pi$, the Lemma \ref{lemmaM} yields there is $C_3>0$ verifying   
$$
\int_{\mathbb{R}^2}\tilde \Phi(\rho_n)dx \leq C_3, \quad \forall n \in \mathbb{N},
$$
showing that $\{\rho_n\}$ is a bounded sequence in $L_{\tilde \Phi}(\mathbb{R}^2)$. From this, the sequence of functionals $\{\tilde \rho_n\}\subset \partial \Psi(u_n)\subset (E_\Phi(\mathbb{R}^{2}))^*$  associated with $\{\rho_n\}$ is also bounded in $(E_\Phi(\mathbb{R}^{2}))^*$, and so, there is $\tilde \rho_0\in (E_\Phi(\mathbb{R}^{2}))^*$, such that $\tilde\rho_n \stackrel{*}{\rightharpoonup} \tilde\rho_0$ in $(E_\Phi(\mathbb{R}^{2}))^*$ for some subsequence, that is, 
\begin{equation}\label{convRHO1}
\int_{\mathbb{R}^2}\rho_n vdx=\langle \tilde \rho_n, v\rangle\to \langle \tilde \rho_0, v\rangle=\int_{\mathbb{R}^2}\rho_0 vdx, \,\, \forall v\in E,
\end{equation}
for some $\rho_0\in L_{\tilde\Phi}(\mathbb{R}^2)$.

Now, using the fact that $\{u_n\}$ is also bounded in $E$, there is $u_\epsilon \in E$ such that 
\begin{equation}\label{convRHO2}
u_n\rightharpoonup u_\epsilon \text{ in } E.
\end{equation}
From (\ref{IG})-(\ref{convRHO2})
\begin{equation}\label{IG01}
0=\int_{\mathbb{R}^2}\nabla u_\epsilon\nabla v +V(x)u_\epsilon v dx-\int_{\mathbb{R}^2}\rho_0 vdx-\epsilon \int_{\mathbb{R}^2}hvdx, \ v\in E.
\end{equation}
To conclude the proof that $u_\epsilon$ is a solution of $(P)$, we must prove that 
\begin{description}
\item [$i)$] $\rho_0(x)\in \partial_t F(x,u_\epsilon(x))$ a.e. in  $\mathbb{R}^2$ and
\item [$ii)$] $\left|[u_\epsilon>t_0]\right|>0$.
\end{description}
To prove the $i)$, we must show that $\{u_n\}$ is strongly convergent to $u_\epsilon$ in $E$, because this fact will imply that $\tilde\rho_0\in\partial \Psi (u_0)$. This way, by Theorem \ref{theoInclusao} 
$$
\rho_0(x)\in \partial_t F(x,u_\epsilon(x)) \text{ a.e in } \mathbb{R}^2.
$$

Related to the second item, the proof is as follows:  If $t_0=0$, then $|[u_\epsilon>t_0]|>0$, because $u_\epsilon \geq 0$ and $u_\epsilon \not= 0$. Next, we will consider the case $t_0 \in (0,t_{\ast})$. Once $\rho_0,u_\epsilon\geq 0$, it follows that 
\begin{eqnarray*}
0&=&\parallel u_\epsilon\parallel^2-\int_{\mathbb{R}^2}\rho_0 u_\epsilon dx-\epsilon \int_{\mathbb{R}^2}hu_\epsilon dx\\
&\leq& \parallel u_\epsilon\parallel^2-\epsilon \int_{\mathbb{R}^2}hu_\epsilon dx,
\end{eqnarray*}
that is, 
\begin{equation}\label{eqtp01}
 \parallel u_\epsilon\parallel^2\geq\epsilon \int_{\mathbb{R}^2}hu_\epsilon dx.
\end{equation}
Arguing by contradiction, we assume that $|[u_\epsilon >t_0]|=0$, for some $t_0\in(0,t_{\ast})$. Thereby, 
$$
f(x,u_\epsilon(x))=0, \text{ a.e. in } \mathbb{R}^2,
$$
from where it follows that 
$$
\partial_t F(x,u_\epsilon(x))=\{0\} \text{ a.e. in  }\mathbb{R}^2.
$$
Consequently,
$$
\rho_0(x)=0 \text{ a.e. in }\mathbb{R}^2.
$$
On the other hand, by Lemma \ref{Inf} and (\ref{eqtp01}),
$$
0>-\delta>I_\epsilon(u_\epsilon)=\frac{1}{2}\parallel u_\epsilon \parallel^2-\epsilon\int_{\mathbb{R}^2}h u_\epsilon dx\geq -\frac{1}{2}\epsilon \int_{\mathbb{R}^2}hu_\epsilon dx \geq -\frac{t_{0}}{2}\epsilon \int_{\mathbb{R}^2}hdx\ ,
$$
implying that
$$
t_{0}\geq\frac{2\delta}{\epsilon \int_{\mathbb{R}^2}hdx}=t_{\ast},
$$
which is a contradiction. \\

\noindent {\bf Convergence of $\{u_n\}$ to $u_\epsilon$ in $E$: } \, Hereafter, fix $\gamma_n:=u_n-u_0$ and recall that $\gamma_n\rightharpoonup 0$ in $E$. By a direct computation,
$$
\parallel u_n\parallel^2 = \parallel u_\epsilon\parallel^2 + \parallel \gamma_n\parallel^2 +o_n(1).
$$
Moreover, we also have 
\begin{eqnarray}\label{eqwn2}
o_n(1)=\langle w_n,u_n\rangle &=& \parallel u_n\parallel^2-\int_{\mathbb{R}^2}\rho_n u_n dx -\epsilon\int_{\mathbb{R}^2} h u_n dx\nonumber\\
&&-\parallel u_\epsilon \parallel^2+\int_{\mathbb{R}^2}\rho_0 u_\epsilon dx +\epsilon\int_{\mathbb{R}^2} h u_\epsilon dx\nonumber\\
&=&\parallel \gamma_n\parallel^2+\left(\int_{\mathbb{R}^2}\rho_0u_\epsilon dx
-\int_{\mathbb{R}^2}\rho_n u_\epsilon dx\right)\nonumber\\
&&+\left(\int_{\mathbb{R}^2}\rho_n u_\epsilon dx-\int_{\mathbb{R}^2}\rho_n u_ndx\right) +o_n(1)\nonumber\\
&=&\parallel \gamma_n\parallel^2-\int_{\mathbb{R}^2}\rho_n (u_n-u_\epsilon)dx+o_n(1)\nonumber\\
&=&\parallel \gamma_n\parallel^2-\int_{\mathbb{R}^2}\rho_n \gamma_ndx+o_n(1).
\end{eqnarray}
On the other hand, by $(f_1)$,
\begin{eqnarray*}
\left|\int_{\mathbb{R}^2}\rho_n \gamma_ndx \right|&\leq& c_1\int_{\mathbb{R}^2}|u_n||\gamma_n| dx+c_2\int_{\mathbb{R}^2}|\gamma_n|\left(e^{\alpha|u_n|^2}-1\right)dx\nonumber\\
&\leq& c_1|u_n|^2_2|\gamma_n|_2^2+{c_2}\int_{\mathbb{R}^2}|\gamma_n|\left(e^{\alpha|\gamma_n|^2}-1\right)dx.\nonumber\\
\end{eqnarray*}
Once $\alpha r^{2} < 4 \pi$, there is $q>1$ close to 1, such that 
$$
M=\sup_{n \in \mathbb{N}}\left(\int_{\mathbb{R}^{2}}\left(e^{\alpha|\gamma_n|^2}-1\right)^{q}dx\right)^{\frac{1}{q}}<+\infty.
$$
Thus, by Lemma \ref{lemmaM} and H\"older inequality 
$$
\left|\int_{\mathbb{R}^2}\rho_n \gamma_ndx\right|\leq c_1|u_n|^2_2|\gamma_n|_2^2+C|\gamma_n|_{q'},
$$
where $\frac{1}{q}+\frac{1}{q'}=1$. Since the embeddings $E \hookrightarrow L^2(\mathbb{R}^2)$ and $E \hookrightarrow L^{q'}(\mathbb{R}^2)$ are compact, we can ensure that
\begin{equation}\label{eqwn3}
\int_{\mathbb{R}^2}\rho_n \gamma_ndx\to 0.
\end{equation}
From (\ref{eqwn2}) and (\ref{eqwn3}), $\gamma_n \to 0$ in 
$E$, or equivalently $u_n\to u_\epsilon$ in $E$, finishing the proof. 

 $\hfill{\rule{2mm}{2mm}}$
 \section{Existence of solution via Mountain Pass}
In this section, we will assume more some conditions on function $f$, namely  $(f_0), (f_2)-(f_6)$. 
By $(f_3)$, there are $\tilde\epsilon,\tilde\delta:=\tilde\delta_{\tilde \epsilon}>0$, satisfying 
$$
2\max\left\{|\xi|; \xi\in\partial_t F(x,t)\right\}<(\lambda_1-\tilde \epsilon)|t|, \ \text{for} \, |t|\leq \tilde\delta \text{ and }x\in \mathbb{R}^2.
$$
From Lebourg's Theorem, there are $\theta(t)\in[0,t]$, with $|t|\leq \tilde\delta$ and $\xi_0\in\partial F_t(x,\theta)$ verifying 
\begin{equation}\label{eqf201}
|F(x,t)|=|F(x,t)-F(x,0)|=|\xi_0||t-0|\leq (\lambda_1-\tilde \epsilon)|t|, \ x\in \mathbb{R}^2.
\end{equation}
Now, by $(f_0)$, given  $q\geq 2$ and $\alpha > \alpha_0$, there is $C=C(q,\tilde\delta)>0$ such that 
$$
|\xi|\leq C|t|^{(q-1)}\left(e^{\alpha |t|^2}-1\right), \ \xi\in \partial_t F(x,t), \ |t|\geq \tilde \delta \text{ and }x\in\mathbb{R}^2.
$$
Applying again Lebourg's Theorem
$$
|F(x,t)|\leq C|t|^q\left(e^{\alpha |t|^2}-1\right), \ \forall t\in\mathbb{R} \text{ and } \forall x\in \mathbb{R}^2.
$$
From this, for $u\in E$ with $u\neq 0$ and $\parallel u\parallel:=\eta_0<\sqrt{\frac{4\pi}{\alpha_0}}$, we see that 
\begin{equation}\label{eqf202}
\int_{\mathbb{R}^2}|F(x,u)|dx\leq C(q,\tilde\delta,\alpha_0)\parallel u\parallel^q.
\end{equation}
Here, we have fixed $\alpha$ close to $\alpha_0$ of the a way that $\alpha \eta_0^{2} < 4 \pi$.

\begin{lemma}\label{lemGMP02}
Assume that $(f_0)$ and $(f_2)-(f_6)$ hold. Then, there exists $\varphi_0\in B^c_{r}(0)$ such that
$$
I_\epsilon(\varphi_0)<\inf_{\parallel u\parallel=r}I_\epsilon(u), \ \epsilon\in\left (0,\epsilon_0\right],
$$
where  $r$ and $\epsilon_0$ are given in Lemma \ref{Inf}.
\end{lemma}
{\bf Proof.} \, Let $\psi_0\in C_0^{\infty}(\mathbb{R}^2)\setminus\{0\}$, $\psi_0>0$, with $supt(\psi_0)\subset K$, where $K\subset\mathbb{R}^2$ is the compact set fixed in $(f_4)$. In this case, for any $\epsilon>0$, 
$$
I_\epsilon(t \psi_0)\leq \frac{t^2}{2}\parallel \psi_0\parallel^2-c_3t^\nu\int_{\mathbb{R}^2}\psi_0^\nu+c_4|K|-t\epsilon\int_{\mathbb{R}^2} h\psi_0 dx,
$$
from where it follows that
$$
\lim_{t \to +\infty}I_\epsilon(t \psi_0)=-\infty.
$$
Thus, the lemma follows choosing $\varphi_0:=t\psi_0\in B^c_{r}(0)$ with $t$ large enough. $\hfill{\rule{2mm}{2mm}}$

\vspace{0.5 cm}

From Lemmas \ref{Inf} and \ref{lemGMP02}, we can use the Mountain Pass Theorem to get a sequence $\{v_n\}\subset E$ verifying 
\begin{equation}\label{seqPS}
I_\epsilon(v_n)\to d_\epsilon \text{ in }\mathbb{R} \text{ and }\lambda_\epsilon(v_n):=\max\{\parallel\xi\parallel_*/ \xi\in\partial I_\epsilon(v_n)\}\to 0,
\end{equation}
where
$$
d_\epsilon:=\displaystyle\inf_{\gamma\in\Gamma}\max_{t\in [0,1]}I_\epsilon(\gamma(t))\quad (\mbox{mountain pass level})
$$
and
$$
\Gamma:=\{\gamma\in C\left([0,1];E\right)/\gamma(0)=0 \text{ and }\gamma(1)=\varphi_0\}.
$$

In the sequel. we intend to show that $I_\epsilon$ verifies the $(PS)_{d_\epsilon}$ condition if the parameter $\mu$ given in $(f_6)$ is large enough. To this end, we need of the following lemma
\begin{lemma}\label{LemaPS1}
Let $\{v_n\}$ be the sequence obtained in (\ref{seqPS}). Then, $\{v_n\}$ is bounded in $E$ and 
$$
\displaystyle\limsup_{n\to\infty}\parallel v_n\parallel\leq \frac{\left(\frac{\tau-1}{\tau}\right)\epsilon+\sqrt{\epsilon^2\left(\frac{\tau-1}{\tau}\right)^2+2d_\epsilon\left(\frac{\tau-2}{\tau}\right)}}{\left(\frac{\tau-2}{\tau}\right)},
$$
where $\tau$ is given in $(f_5)$.
\end{lemma}
{\bf Proof.} \, Let $w_n \in E^{*}$ and $\rho_n\in \partial \Psi(v_n) $ verifying  
$$
\parallel w_n\parallel_*=\lambda_\epsilon(v_n) \quad \mbox{and} \quad  \langle w_n,v_n\rangle=\parallel v_n\parallel^2 -\int_{\mathbb{R}^2}\rho_n v_ndx-\epsilon\int_{\mathbb{R}^2}h v_ndx.
$$
From  $(f_5)$ and Theorem \ref{theoInclusao}, 
\begin{eqnarray}\label{eqPS01}
d_\epsilon+o_n(1)+o_n(1)\parallel v_n\parallel&\geq &I_\epsilon(v_n)-\frac{1}{\tau}\langle w_n,u_n\rangle\nonumber\\
&=&\left(\frac{1}{2}-\frac{1}{\tau}\right)\parallel v_n\parallel^2+\int_{\mathbb{R}^2}\frac{1}{\tau}\rho_nv_n-F(x,v_n)dx\nonumber\\
&&+\left(\frac{1}{\tau}-1\right)\epsilon\int_{\mathbb{R}^2}h v_ndx\\
&\geq&\left(\frac{1}{2}-\frac{1}{\tau}\right)\parallel v_n\parallel^2+\left(\frac{1}{\tau}-1\right)\epsilon \|h\|_* \|v_n\|,\nonumber
\end{eqnarray}
which implies that $\{v_n\}$ is bounded in $E$. Moreover, as $\{v_n\}$ does not converge to $v=0$ in $E$, we can assume that for some subsequence, 
$$
 l:=\displaystyle\lim_{n\to\infty}\parallel v_n\parallel>0.
$$
Consequently, by  (\ref{eqPS01}),
$$
d_\epsilon+\left(\frac{\tau-1}{\tau}\right)\epsilon\parallel h\parallel_*l\geq  \left(\frac{1}{2}-\frac{1}{\tau}\right) l^2,
$$
that is
$$
\left(\frac{1}{2}-\frac{1}{\tau}\right) l^2-\left(\frac{\tau-1}{\tau}\right)\epsilon\parallel h\parallel_*l-d_\epsilon \leq 0.
$$
As $l>0$, we must have 
$$
l\leq \frac{\left(\frac{\tau-1}{\tau}\right)\epsilon+\sqrt{\epsilon^2\left(\frac{\tau-1}{\tau}\right)^2+4d_\epsilon\left(\frac{\tau-2}{2\tau}\right)}}{2\left(\frac{\tau-2}{2\tau}\right)},
$$
which completes the proof. $\hfill{\rule{2mm}{2mm}}$
\begin{lemma}\label{LemaPS2}
Assume $(f_0)-(f_6)$. Then, there are $\epsilon_1, \mu^*>0$ and $t_1>0$ such that
$$
\frac{\left(\frac{\tau-1}{\tau}\right)\epsilon+\sqrt{\epsilon^2\left(\frac{\tau-1}{\tau}\right)^2+2d_\epsilon\left(\frac{\tau-2}{\tau}\right)}}{\left(\frac{\tau-2}{\tau}\right)} < \sqrt{\frac{4 \pi}{\alpha_0}}.
$$
for all $\epsilon \in (0, \epsilon_1), \mu \geq \mu^*$ and $t_0 \in [0,t_1)$. 
\end{lemma}
{\bf Proof.} \, Consider the function $\psi_0$ used in the proof of Lemma \ref{lemGMP02}. Then, 
$$
\sup_{t \in [0,t_0]}I_\epsilon(t \psi_0)\leq \frac{t_0^{2}}{2}\|\psi_0\|^{2},
$$
and so, there is $t_2>0$ such that
$$
\sup_{t \in [0,t_0]}I_\epsilon(t \psi_0)\leq \epsilon^{2}
$$
for $t_0 \in [0, t_2)$. On the other hand, by $(f_6)$,
$$
\sup_{t \geq t_0}I_\epsilon(t \psi)\leq \displaystyle \max_{t\geq 0}\left\{\frac{t^2}{2}\parallel \psi_0\parallel^2-\mu t^{p}\int_{\mathbb{R}^2}\psi_0^pdx\right\}+\mu c_2t_0^{p}|supt(\psi_0)|, 
$$
that is,
$$
\sup_{t \geq t_0}I_\epsilon(t \psi_0)\leq \left(\frac{1}{2 p^{\frac{2}{p-2}}}-\frac{1}{p^{\frac{p}{p-2}}}\right)\frac{1}{\mu^{\frac{2}{p-2}}}\left(\frac{\parallel \psi_0\parallel}{|\psi_0|_p}\right)^{\frac{2p}{p-2}} + \mu c_2t_0^{p}|supt(\psi_0)|.
$$
Now, fix $\mu^*>0$ such that 
$$
\left(\frac{1}{2 p^{\frac{2}{p-2}}}-\frac{1}{p^{\frac{p}{p-2}}}\right)\frac{1}{\mu^{\frac{2}{p-2}}}\left(\frac{\parallel \psi_0\parallel}{|\psi_0|_p}\right)^{\frac{2p}{p-2}} \leq \epsilon^{2}, \quad \forall \mu \geq \mu^*
$$
and $t_3=t_3(\mu,\epsilon)>0$ such that
$$
\mu c_2 t_0^{p}|supt(\psi_0)| \leq \epsilon^{2}, \quad \forall t \in [0,t_3].
$$ 
From this, for  $t_1=\min\{t_2,t_3\}$, we must have  
$$
\sup_{t \geq t_0}I_\epsilon(t \psi_0)\leq 2\epsilon^{2},
$$ 
and so,
$$
d_\epsilon \leq \max_{t \geq 0}I_\epsilon(t\psi_0)\leq 2 \epsilon^{2}. 
$$
Hence, there is $c_1>0$ independent of $\epsilon$ such that 
$$
\frac{\left(\frac{\tau-1}{\tau}\right)\epsilon+\sqrt{\epsilon^2\left(\frac{\tau-1}{\tau}\right)^2+2d_\epsilon\left(\frac{\tau-2}{\tau}\right)}}{\left(\frac{\tau-2}{\tau}\right)} \leq c_1 \epsilon.
$$
Then, there is $\epsilon_0>0$ such that
$$
\frac{\left(\frac{\tau-1}{\tau}\right)\epsilon+\sqrt{\epsilon^2\left(\frac{\tau-1}{\tau}\right)^2+2d_\epsilon\left(\frac{\tau-2}{\tau}\right)}}{\left(\frac{\tau-2}{\tau}\right)} < \sqrt{\frac{4 \pi}{\alpha_0}}, \quad \forall \epsilon \in (0,\epsilon_0).
$$
$\hfill{\rule{2mm}{2mm}}$

As an immediate consequence of the last lemma, we have the following corollary.

\begin{coro} \label{convergencia} Let $\{v_n\}$ be the sequence obtained in (\ref{seqPS}). Then, there is $\epsilon_0$ such that
$$
\limsup_{n \to +\infty}\|v_n\|^{2}< \frac{4 \pi}{\alpha_0}, \quad \forall \epsilon \in (0, \epsilon_0).
$$
Moreover, there is a subsequence of $\{v_n\}$ still denoted by itself, and $v_\epsilon \in E$ such that $v_n \to v_\epsilon$ in $E$.
\end{coro}
\noindent {\bf Proof.} \, The first part of the lemma is an immediate consequence of Lemmas \ref{LemaPS1} and \ref{LemaPS2}. The proof of the second part follows the same idea explored in the proof of  Theorem \ref{SolEkland}. $\hfill{\rule{2mm}{2mm}}$

\begin{theorem}\label{mountainpass}
Assume $(V_1)-(V_2)$ and $(f_0)-(f_6)$. Then, there are $\epsilon_0, \mu^*$ and $t_1>0$, such that problem $(P)$ possesses a solution $v_\epsilon \in E $, with $I_\epsilon(v_\epsilon)=d_\epsilon >0$, for all $\epsilon \in (0,\epsilon_0),$ $t_0\in [0,t_{1})$ and $\mu \geq \mu^*$. Moreover, decreasing $\epsilon_0$ and $t_1$, and   increasing $\mu^*$, if necessary, we have two solutions $u_\epsilon, v_\epsilon \in E$ with
$$
I_\epsilon(u_\epsilon)=c_\epsilon<0<d_\epsilon=I_\epsilon(v_\epsilon).
$$ 
\end{theorem}
\noindent {\bf Proof.}\, The theorem follows applying the Lemmas \ref{Inf} and \ref{lemGMP02} and Corollary \ref{convergencia}. $\hfill{\rule{2mm}{2mm}}$

\vspace{1 cm}

\noindent Claudianor O. Alves and Jefferson A. Santos \\
\noindent Universidade Federal de Campina Grande, \\
\noindent Unidade Acad\^emica de Matem\'atica ,\\
\noindent CEP:58109-970, Campina Grande - PB, Brazil\\
\noindent e-mail: coalves@mat.ufcg.edu.br and jefferson@mat.ufcg.edu.br \\

\end{document}